\documentclass[journal]{IEEEtran}
\usepackage{graphicx}
\usepackage{bm}
\usepackage{bbm}
\usepackage{caption}
\usepackage{amsmath}
\usepackage{amssymb}
\usepackage{mathtools}
\usepackage{multirow}
\usepackage{subcaption}
\usepackage{yfonts}
\usepackage{xcolor}
\usepackage[ruled]{algorithm2e}
\usepackage[margin=0.75in,bottom=0.8in]{geometry}
\usepackage{titling}
\date{}
\newtheorem{remark}{Remark}
\newtheorem{prob}{Problem}
\newtheorem{property}{Property}

\title{Solutions to Differential Algebraic Inequalities with Composite Bernstein Polynomials}
\author{Maxwell Hammond, Gage MacLin, Laurent Jay, Venanzio Cichella}

\begin{document}
\maketitle

\begin{abstract}
The Bernstein polynomial basis sees significant use owing to its unique properties, particularly in the field of optimal control. However, the basis is known to have a slow rate of convergence to the function it approximates. With this in mind, we introduce two collocation methods for solving general ordinary differential equations using composite Bernstein polynomials to preserve the basis properties while improving convergence. Of particular note is the integration based method which uses a minimal number of variables to describe the resulting composite polynomial, reducing computational effort. In addition, we exploit the convex hull property of the Bernstein polynomial basis in order to enforce inequality constraints in differential algebraic inequalities, highlighting the benefits of the basis in function approximation. Solutions to six numerical examples are provided as well as discussion of the advantages and disadvantages of the proposed solution methodologies.
\end{abstract}

\section{Introduction}
In spite of their relatively slow convergence rate for approximating functions, Bernstein polynomials have worked their way into a number of fields owing to their numerical stability and desirable properties \cite{farouki2012bernstein}.  In particular, the convex hull property of this basis has led to its use in the field of optimal control, where inequality constraints can be safely enforced through the polynomial coefficients, even at low orders of approximation \cite{cichella2017optimal,kielas2022bernstein}. Research in this field may be extended by the use of efficient numerical methods to solve complex system dynamics which take the form of differential equations \cite{rao2009survey}. This, and other similar cases have motivated efforts to deploy Galerkin methods \cite{bhatti2007solutions,doha2011derivatives}, Adomian decomposition \cite{qasim2018adomian}, collocation \cite{dascioglu2013bernstein}, and other numerical methods using the Bernstein basis to approximate the solutions of linear and nonlinear differential equations. Nevertheless, these methods struggle to compare with the convergence rate and speed of well established methods such as Runge-Kutta methods \cite{wanner1993solving,wanner1996solving}.

While the use of Bernstein basis may not achieve the pinnacle of efficiency for modern approaches to solving ordinary differential equations (ODEs), improvements can be made to approaches taken with this basis and their applications can be expanded. While the use of composite polynomials is not foreign within numerical methods, especially for finite element cases  \cite{jacobi2010minimization,tadmor2012review}, minimal attention has been given to composite Bernstein polynomials (CBPs) in literature surrounding differential equations. Instead, focus is placed on geometry and data fitting applications \cite{jacobi2010minimization,gousenbourger2019data}, with some notable exceptions in optimal path generation \cite{mehdi2015avoiding,cao2022well,patterson2023hermite,maclin2024optimal}.  Broadening slightly the scope to highlight the advantageous properties of the Bernstein basis, differential algebraic inequalities (DAIs) are a class of problem made up of a set of differential algebraic equations (DAEs) \cite{peter2010alpha,spiteri2000programming}. These problems often appear in engineering contexts describing dynamic systems, and solution methodologies would benefit both from increased efficiency in solving differential equations and enforcing inequalities. Within this manuscript we seek to exploit the particular properties of CBPs in order to rewrite ordinary differential equations (ODEs) as sets of discrete algebraic equations which can be solved with widely available software. To this end, we present two collocation methods which provide solutions to ODEs in the form of CBPs.  Further, we present methods for solving DAIs within the same mathematical framework. This is accomplished by exploiting the convex hull property of Bernstein polynomials, leading to a robust and safe framework for satisfying these inequalities.

The paper is organized as follows: Section \ref{sect:bern} introduces the Bernstein polynomial basis and its key properties, Section \ref{sect:meth} introduces a derivative and an integral based collocation method using composite Bernstein polynomials to solve ODEs, and introduces a method for enforcement of inequality constraints with CBPs, Section \ref{sect:num} gives some numerical examples of ODEs and DAIs, and finally Section \ref{sect:conc} provides some conclusion.

\section{The Bernstein Basis} \label{sect:bern}
The Bernstein basis of $n$-th degree over the interval $[s_0,s_f]$ is given by
\begin{equation} \label{eq:basis}
    B_{i}^{n}(s)={n \choose i}\frac{\left(s-s_0\right)^i\left(s_f-s\right)^{n-i}}{\left(s_f-s_0\right)^n}, \quad i = 0,1,2,\dots n,
\end{equation}
where $${n \choose i} = \frac{n!}{i!\left(n-i\right)!}.$$ An $n$-th order univariate polynomial in Bernstein form is expressed 
\begin{equation} \label{eq:bpoly}
    x_n(s) = \sum_{i=0}^{n} \bar{x}_i^nB_i^n(s), \quad s \in \left[s_0,s_f\right].
\end{equation}
where $\bar{x}^n\in\mathbb{R}^{n+1}$ is the vector of polynomial coefficients, or control points (CP). For brevity, we refer to \eqref{eq:bpoly} as a Bernstein polynomial (BP).

The Bernstein basis has many relevant properties and what follows is an introduction of those of particular importance to this manuscript. A more complete list of the properties can be found in \cite{farouki2012bernstein}. 
\vspace{6pt}
\begin{property}[non-negativity and unity] \label{prop:unity}
    The basis functions \eqref{eq:basis} satisfy $B_i^n(s) \ge 0$ for $s \in [s_0,s_f],$ and $\sum_{i=0}^{n} B_i^n(s) = 1$.
\end{property}
\vspace{6pt}
\begin{property}[end-point values] \label{prop:end}
    The basis functions \eqref{eq:basis} satisfy $B_0^n(s_0) = B_n^n(s_f) = 1$. When considering the BP \eqref{eq:bpoly} and Property \ref{prop:unity}, this leads to the relations
    \begin{equation*}
        x_n(s_0) = \bar{x}_0^n, \quad x_n(s_f) = \bar{x}_n^n.
    \end{equation*}
\end{property}
\vspace{6pt}
\begin{property}[Arithmetic operations] \label{prop:arithm}
Addition and subtraction between two Bernstein polynomials or surfaces can be performed directly through the addition and subtraction of their control points. 
The control points of the Bernstein polynomial \(x_{m+n}(s)\) resulting from the multiplication two Bernstein polynomial, \(p_m(s)\) and \(q_n(s)\) with control points $\bar{p}^{m}$ and $\bar{q}^{n}$ can be obtained by 
\begin{equation}
    \bar{x}_k^{m+n}=\sum_{j=max(0,k-n)}^{min(m,k)}\frac{{m \choose j}{n \choose k-j}}{{m+n \choose k}} \bar{p}_j^m\bar{q}_{k-j}^n
\end{equation}
\end{property}
\vspace{6pt}
\begin{property}[derivatives] \label{prop:deriv}
    The derivative of the basis function \eqref{eq:basis} is
    \begin{equation} \label{eq:deriv}
        \frac{d B_{i}^{n}(s)}{ds} = \frac{n}{s_f-s_0}(B_{i-1}^{n-1}(s)-B_{i}^{n-1}(s))
    \end{equation}
    where $B_{-1}^{n-1}(s)=B_{n}^{n-1}(s)=0$. Thus, the derivative of \eqref{eq:bpoly} can be obtained
    \begin{equation*}
        \begin{split}
            x'_{n-1}(s) = \sum_{i=0}^{n-1} \bar{x}_i^{',n-1} B_i^{n-1}(s), \\ \bar{x}^{',n-1} = \bar{x}^n\boldsymbol{\Delta}_n,
        \end{split}
    \end{equation*}  
    where 
    \begin{equation*} 
        \boldsymbol{\Delta}_n= \frac{n}{s_f-s_0}
        \begin{bmatrix}
            -1 & 0 & \ldots & 0 \\
            1 & \ddots & \ddots & \vdots \\
            0 & \ddots & \ddots & 0 \\
            \vdots & \ddots & \ddots & -1 \\
            0 & \ldots & 0 & 1 
        \end{bmatrix} \in \mathbb{R}^{(n+1)\times n}.
    \end{equation*}
\end{property}
\vspace{6pt}
\begin{property}[integrals] \label{prop:int}
        The integral of the basis function \eqref{eq:basis} is
    \begin{equation*}
        \int B_i^n(s) ds = \frac{s_f-s_0}{n+1}\sum_{k=i+1}^{n+1}B_k^{n+1}(s), 
    \end{equation*}
    leading to the indefinite integral of the function \eqref{eq:bpoly} being given by 
    \begin{equation*}
        \begin{split}
             \int x_n'(s) ds = \sum_{i=0}^{n+1}\bar{x}^{n+1}_iB_i^{n+1}(s), \\
             \bar{x}^{n+1} = \bar{x}^{',n}\boldsymbol{\gamma}_n + x_n(s_0) \mathbbm{1}_{n+2}^\top,
        \end{split}
    \end{equation*}
    where $x_n(s_0)$ is an integration constant and $\mathbbm{1}_{n+2}\in \mathbb{R}^{n+2}$ is a vector of ones 
    and 
    \begin{equation} \label{eq:single_int}
        \boldsymbol{\gamma}_n = \frac{s_f-s_0}{n+1}\begin{bmatrix}
        0 & 1 & \dots & 1 \\
         \vdots  & \ddots & \ddots &\vdots \\
        0  & \dots  & 0 & 1
        \end{bmatrix}  \in \mathbb{R}^{n+1\times n+2}.
    \end{equation}
\end{property}
\vspace{6pt}
\begin{property}[Degree Elevation] \label{prop:el}
    For ${n_e}>n$, a Bernstein polynomial of degree $n$ with control points $\bar{x}^n$ can be degree elevated to order $n_e$ by
    \begin{equation*}
        \bar{x}^{n_e}=\bar{x}^n\textbf{E}_n^{n_e}
    \end{equation*}
    where $\textbf{E}_n^{n_e}=\{e_{j,k}\}\in\mathbb{R}^{(n+1)\times({n_e}+1)}$
    \begin{equation*} \label{eq:degel}
    \begin{split}
            e_{i,i+j}= \left\{ \begin{matrix}
                \frac{{{n_e}-n \choose j}{n \choose i}}{{{n_e} \choose i+j}} & j\le {n_e}-n, \ i\le n \\ \\ 
                0 & \text{otherwise}
            \end{matrix} \right. 
    \end{split}
    \end{equation*}
\end{property}
\vspace{6pt}
\begin{remark} \label{rem:order}
    It is clear by Property \ref{prop:deriv} that the operations for obtaining the derivative of a Bernstein polynomial reduces the degree of that polynomial by one. Similarly, by Property \ref{prop:int}, the integration operation increases the polynomial order by one. In the derivative case, Property \ref{prop:el} can be used to generate a derivative matrix which does not reduce the polynomial order 
    \begin{equation}
        \begin{split}
            \bar{x}_i^{',n} = \bar{x}^n\textbf{D}_n, \\
            \textbf{D}_n = \boldsymbol{\Delta}_n \textbf{E}_{n-1}^n.
        \end{split}
    \end{equation}
\end{remark}
\vspace{6pt}
\begin{property}[Convex Hull] \label{prop:convex}
A Bernstein polynomial lies within the convex hull defined by its control points 
\begin{equation*}
    \min_{i\in\{0,\dots,n\}}\bar{x}_i^n \le x_n(s) \le \max_{i\in\{0,\dots,n\}}\bar{x}_i^n. 
\end{equation*}
\end{property}

\section{Methods}  \label{sect:meth}
\begin{prob} \label{prob:gen}
    Consider the general $r$th-order nonlinear ODE
    \begin{equation} \label{eq:gen_ode}
        \begin{split}
             x^{(r)}(s)=f(s,x(s),x'(s),\dots, x^{(r-1)}(s)), \\ s_0 \le s \le s_f,   
        \end{split}
    \end{equation}
    with initial conditions 
    \begin{equation} 
        \begin{split}
            \sum_{l=0}^{r-1}\lambda_{jl}x^{(l)}(a)=\mu_j,\\
            j=0,1,\dots,r-1, \ a\in[s_0,s_f],
        \end{split}
    \end{equation}
    or boundary conditions
    \begin{equation} \label{eq:BC}
        \begin{split}
            \sum_{l=0}^{r-1}[\alpha_{jl}x^{(l)}(s_0)+\beta_{jl}x^{(l)}(s_f)]=\gamma_j,\\
            j=0,1,\dots,r-1,
        \end{split}
    \end{equation}
    where $\lambda_{jl}, \ \mu_{j}, \ \alpha_{jl}, \ \beta_{jl},$ and $\gamma_j$ are known constants, and $x(s)$ is an unknown function. 
\end{prob}

In what follows, we present two methods which discretize Problem \ref{prob:gen} into solvable sets of algebraic equations. In both cases, the approximation of the solution, $x(s)$, takes the form of a set of $K$ composite Bernstein polynomials (CBPs) over the sub-intervals $[s_i,s_{i+1}], \ i=0,\dots,K-1$ where $s_0<s_1<\dots<s_{K}$ and $s_{K}=s_f$, i.e.
\begin{equation} \label{eq:comp}
    \begin{split}
        X(s) = \begin{cases}
            x_n^{[0]}(s)=\sum_{i=0}^{n} \bar{x}_i^{n[0]} B_i^n(s), \  s\in[s_0,s_1]\\
            x_n^{[1]}(s)=\sum_{i=0}^{n} \bar{x}_i^{n[1]} B_i^n(s), \ s\in[s_1,s_2]\\
            \vdots \\
            \begin{split}
                x_n^{[K-1]}(s)=\sum_{i=0}^{n} \bar{x}_i^{n[K-1]} &B_i^n(s), \\  &s\in[s_{K-1},s_{f}]\\    
            \end{split}
        \end{cases},
    \end{split}
\end{equation}
which is reminiscent of a splines. For ease of manipulation, the full set of control points (CP) of the CBPs will be written $\bar{X}^{n,K}=[\bar{x}^{n[0]},\bar{x}^{n[1]},...,\bar{x}^{n[k-1]}]\in\mathbb{R}^{K(n+1)}$. Here, bracketed numbers in the superscripts of control points denote which polynomial in a CBP is being referred to. For simplicity, each Bernstein polynomial component of the composite Bernstein polynomial is considered to be of the same order, \(n\). However, the methods described can be readily extended to cases where the polynomials are of different orders.

\subsection{Control Point Collocation} \label{sect:cpcoll}
\begin{figure}
    \centering
    \includegraphics[width=0.9\linewidth]{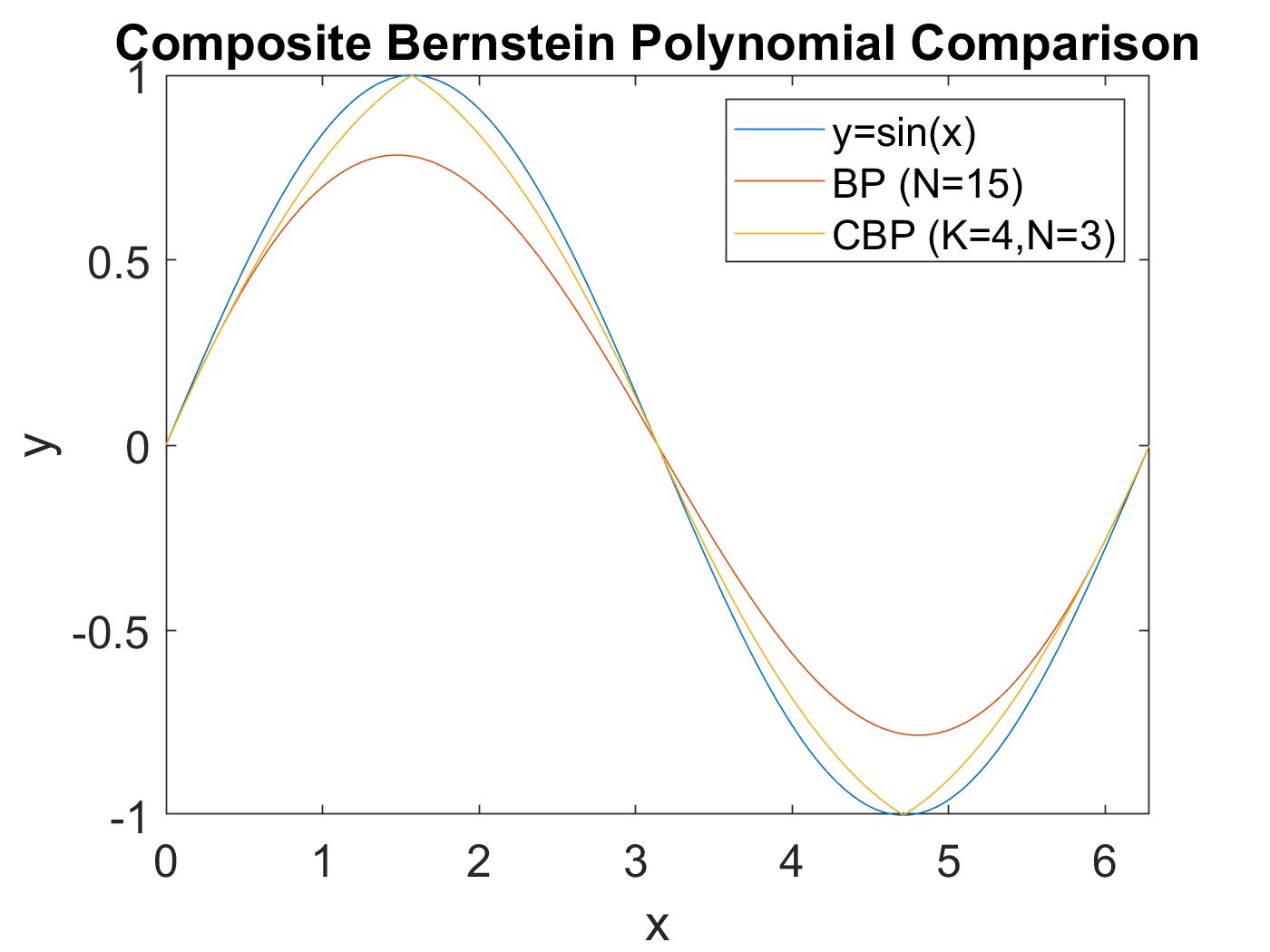}
    \caption{Approximation of $y=\sin(x)$ with BP and CBP with 16 control points. Control points are directly sampled from the function.}
    \label{fig:bp_cbp_comp}
\end{figure}

Traditionally, the Bernstein approximation of a function is performed by choosing control points (CPs) that directly sample the function at equidistant nodes. This approach was originally conceived by Sergei Bernstein when he developed the polynomials to provide a constructive proof of the Weierstrass approximation theorem \cite{bernstein1912demonstration}. This is also the preferred method over the use of the Bernstein basis matrix owing to its stability as the order of the polynomial $n$ increases. A similar practice can be adopted in a composite case; however, CBPs converge to the approximated function with order $\textbf{O}(\frac{1}{nK^2})$ where a single BP converges with $\textbf{O}(\frac{1}{n})$. This improved convergence is illustrated by the approximation of a simple function in Figure \ref{fig:bp_cbp_comp}. 
With this in mind, we employ a method which directly collocates the control points $\bar{X}^{n,K}$ on the ODE in Problem \ref{prob:gen} over each sub-interval $[s_i,s_{i+1}]$. 

Extending Properties \ref{prop:deriv} and \ref{prop:el} (see Remark \ref{rem:order}) for the CBP case, the CPs of the first $r$ derivatives of $X(s)$ with equivalent order can be generated with the block diagonal derivative matrix
\begin{equation}
    \begin{split}
        \bar{X}^{',n,K} = \bar{X}^{n,K}\textbf{D}_{n,K}, \\
        \textbf{D}_{n,K} = \begin{bmatrix}
            \textbf{D}_n^{[0]} & 0 & \dots & 0 \\
            0 & \ddots & \ddots & \vdots \\
            \vdots & \ddots & \ddots & 0 \\
            0 & \dots & 0 & \textbf{D}_n^{[K-1]}
        \end{bmatrix},    
    \end{split}
\end{equation}
where $\textbf{D}_{n,K}\in\mathbb{R}^{K(n+1)\times K(n+1)}$ and $\textbf{D}_n^{[i]}$ is the differentiation matrix associated with the sub-interval $s\in[s_i,s_{i+1}]$. These CPs can be used to represent the problem discretely, treating them as points on the function at equidistant nodes in their respective sub-intervals. 
These points will be denoted $s_{n,K}$ for the entire interval $[s_0,s_f]$. Continuity of order up to $C^{r-1}$ must be enforced at the intersections of the composite polynomials (knots) to create sufficient conditions to solve the ODE on each sub-interval. This can be accomplished with equality constraints exploiting Property \ref{prop:end} at the knots of the derivatives of $X(s)$. With the decision variable being the control points $\bar{X}^{n,K}$, the discretized problem can be formulated: 
\begin{prob} \label{prob:deriv}
    \begin{equation} \label{eq:disc_ode}
        \begin{split}
             \bar{X}^{n,K}(\textbf{D}_{n,K})^r &= \\
             f(s_{n,K},\bar{X}^{n,K},&\bar{X}^{n,K}\textbf{D}_{n,K},\dots, \bar{X}^{n,K}(\textbf{D}_{n,K})^{r-1}),  \\
        \end{split}
    \end{equation}
    Subject to initial conditions 
    \begin{equation} \label{eq:disc_IC}
        \begin{split}
            \sum_{l=0}^{r-1}\lambda_{jl}X^{(l)}(a) = \mu_j,  \\
            j = 0,1,\dots,r-1,
        \end{split}
    \end{equation}
    boundary conditions,
    \begin{equation} \label{eq:disc_BC}
        \begin{split}
            \sum_{l=0}^{r-1}[\alpha_{jl}\{\bar{X}^{n,K}(\textbf{D}_{n,K})^l\}_0+\beta_{jl}\{\bar{X}^{n,K}(\textbf{D}_{n,K})^l\}_f]=\gamma_j,\\
            j=0,1,\dots,r-1,
        \end{split}
    \end{equation}
    where $\{\cdot\}_0$ and $\{\cdot\}_f$ represent the first and final CPs of the interior result respectively. The continuity conditions are written
    \begin{equation}\label{eq:ode_cont}
        \begin{split}
            \bar{x}^{n[j]}_n&=\bar{x}^{n[j+1]}_0, \\
            \{\bar{x}^{n[j]}\textbf{D}_n\}_f&=\{\bar{x}^{n[j+1]}\textbf{D}_n\}_0,\\
            &\vdots\\
            \{\bar{x}^{n[j]}(\textbf{D}_n)^{r-1}\}_f&=\{\bar{x}^{n[j+1]}(\textbf{D}_n)^{r-1}\}_0, \\
            j &= 1,..,K-1.
        \end{split}
    \end{equation}
\end{prob}
\begin{remark}
    As a result of the continuity conditions, the order of the polynomial must be selected to satisfy $n > r$. This ensures that the polynomials are of an order greater than zero at the $r-1$ derivative and continuity can thereby be enforced.
\end{remark}
This method leads to a problem formulation with minimal continuity constraints as well as flexibility in the choice of two degrees of freedom, $n$ and $K$. 
However, the number of variables in this case grows according to $K(n+1)$  and the problem can be simplified in cases where this flexibility is unnecessary.

\subsection{Knot Only Collocation} \label{sect:knot_only}
We wish to find an approximate solution to Problem \ref{prob:gen} in the form \eqref{eq:comp} which collocates the differential equation on $K+1$ points, $s_0,s_1,\dots,s_{K}=s_f$, and can be represented by a minimum number of variables. 
Notice that $X(s)$ can be written in terms of its $M-th$ derivative and its initial conditions by
\begin{equation} \label{eq:int_expanded}
    \begin{split}
        X(s) = \int \dots &\int X^{(M)} ds \dots ds \\ 
        &+ X_0 + c_1 X_0' + \dots + c_{M-1} X_0^{(M-1)},
    \end{split}
\end{equation}
where $c_i, \ i=1,\dots, M-1$ are some constants. Using the integration property of Bernstein polynomials, \eqref{eq:int_expanded} can be reformulated as an algebraic operation. Let the $M$-th derivative of the solution, $x^{(M)}(s), \ M = r+1$ be approximated by $K$ zeroth order CBPs 
\begin{equation}
    X^{(M)}(s)=\begin{cases}
        \bar{x}^{(M),0[0]}, \quad s\in[s_0,s_1]\\
        \bar{x}^{(M),0[1]}, \quad s\in[s_1,s_2]\\
        \vdots\\
        \bar{x}^{(M),0[K-1]}, \quad s\in[s_{K-1},s_f]
    \end{cases}.
\end{equation}
The vector of unknowns needed to solve the system is given as 
$$\theta_0=[\bar{X}^{(M),0,K},X_0,X'_0,\dots,X^{(M-1)}_0]$$
where $\theta_0\in\mathbb{R}^{K+M}$ and $X_0,X'_0,\dots,X^{(M-1)}_0$ are the integration constants for the zeroth to $(M-1)$-th derivative of $X(s)$. The integration Property \ref{prop:int} of Bernstein polynomials can be exploited to obtain the control points of the derivatives of $X(s)$ in the form $\theta_n$ by
\begin{equation} \label{eq:theta_int}
    \begin{split}
        \theta_{1}=\theta_0\boldsymbol{\zeta}_{0},  \\
        \vdots \\
        \theta_{M}=\theta_{M-1}\boldsymbol{\zeta}_{M-1},
    \end{split}
\end{equation}
where
\begin{equation*}
        \theta_n= [\bar{X}^{(M-n),n,K},X_0,\dots,X^{(M-1)}_0]       
\end{equation*}
and $\zeta_n$ is a modified integration matrix. Thus, a transformation matrix can be written
\begin{equation} \label{eq:transformation}
    \begin{split}
        \text{knots}(\bar{X}^{(M-m),M-m,K})=\theta_0\textbf{T}_{M-m}, \\
        \textbf{T}_{M-m} = \textbf{P}_{m-1}\boldsymbol{\zeta}_{m-1}\dots\boldsymbol{\zeta}_{0},
    \end{split}
\end{equation}
which gives the $K+1$ knots of the $(M-m)$-th derivative of $X(s)$. The formulation of the matrices $\boldsymbol{\zeta}_n$ and $\textbf{P}_n$ can be found in the Appendix \ref{app:int}. It is important to note that the matrix $\boldsymbol{\zeta}_i$ facilitates the integration of the function from the $i$-th to the $(i-1)$-th derivative of the solution, and simultaneously ensures continuity constraints, i.e., the resulting integrand is continuous.

Problem \ref{prob:gen} can now be reformulated as a collocation problem on the knots of X(s) and its derivatives:
\begin{prob} \label{prob:integ}
\begin{equation}
    \theta_0\textbf{T}_{r} = f(s_i,\theta_0\textbf{T}_{0},\dots,\theta_0\textbf{T}_{r-1}).
\end{equation}
With initial conditions
\begin{equation} 
    \begin{split}
        \sum_{l=0}^{r-1}\lambda_{jl}X^{(l)}(a)=\mu_j,\\
        j=0,1,\dots,r-1, \ a\in[s_0,s_f],
    \end{split}
\end{equation}
and boundary conditions 
\begin{equation}
    \begin{split}
        \sum_{l=0}^{r-1}[\alpha_{jl}\{\theta_0\textbf{T}_l\}_0+\beta_{jl}\{\theta_0\textbf{T}_l\}_f]=\gamma_j,\\
        j=0,1,\dots,r-1,
    \end{split}
\end{equation}
\end{prob}
\begin{remark} \label{rem:add_cond}
    For $M=r+2$, an additional condition is necessary to obtain a solution. One way to obtain this condition is by evaluating the Bernstein polynomial on the first interval $s\in[s_0,s_1]$ and satisfying \eqref{eq:gen_ode} at $\sigma=\frac{s_1-s_0}{2}$, i.e.
    \begin{equation} \label{eq:add_cond}
        X^{(r)}(\sigma)=f(s,X(\sigma),X'(\sigma),\dots, X^{(r-1)}(\sigma)).
    \end{equation}
    Higher order $M$ can be similarly obtained; however further investigation should be conducted on choices of $\sigma$ which yield stable results. 
    Additional consideration could also be given to classical approaches taken in the case of interpolating splines \cite{shikin1995handbook}.
\end{remark}

\subsection{Inequality Constraints} \label{sect:VI}
Consider the general DAI
\begin{prob}
    \begin{equation}
        \begin{split}
            x^{(r)}(s)=f(s,x(s),x'(s),\dots, x^{(r-1)}(s)), \\
            g(s,x(s),x'(s),\dots,x^{(r-1)}(s)) \le 0, \\
            s_0 \le s \le s_f.
        \end{split}
    \end{equation}
    subject to initial and boundary conditions.
\end{prob}
Again, let the solution $x(s)$ be approximated by a CBP in the form \eqref{eq:comp}. Here, the ODE condition is satisfied through the methods described in the above sections. The composite polynomial $X(s)$ is simultaneously checked within the inequality condition leading to 
\begin{equation} \label{eq:comp_ineq}
    \begin{split}
        g(s,X(s),X'(s),\dots,X^{(r)}(s)) \le 0, \\
        s_0 \le s \le s_f.
    \end{split}
\end{equation}
Using the Properties \ref{prop:arithm}, \ref{prop:deriv}, and \ref{prop:el}, the function on the left side of the inequality can be represented as a CBP through the manipulation of the control points $\bar{X}^{n,K}$. Let the resulting CBP have the control points $\bar{G}^{n,K}=[\bar{g}^{n[0]},\bar{g}^{n[1]},...,\bar{g}^{n[k-1]}]$ (similar to \eqref{eq:comp}). By the convex hull Property \ref{prop:convex} of Bernstein polynomials
\begin{equation*}
    \begin{split}
        \min_{i\in\{0,\dots,n\}}\bar{g}_i^{n[j]} \le g_n^{[j]}(s) \le \max_{i\in\{0,\dots,n\}}\bar{g}_i^{n[j]}, \\
        j=0,\dots,K-1.
    \end{split}
\end{equation*}
Thus, the inequality \eqref{eq:comp_ineq} is satisfied when
$\max_{i,j} \bar{G}_{i,j}^{n,K}\le 0$.

\begin{remark}
    The check on the inequality must be performed on all control points to guarantee satisfaction. 
\end{remark}

\section{Numerical Results} \label{sect:num}
In this section we provide four numerical examples for the presented methods of solving ODEs to demonstrate the convergence of each method. Additionally we provide two numerical examples of DAIs using the convex hull property to satisfy inequalities. Solutions are obtained in MATLAB using the fsolve and fmincon functions on an Intel\textsuperscript{\textregistered} Core\textsuperscript{\texttrademark} i9-10885H CPU at 2.40GHz, 2400 Mhz with 8 Cores and 16 Logical Processors. For simplicity in convergence plots, solutions in these examples are generated on equidistant knots, however this is not required for the presented methods. An implementation of these examples and other examples can be found in https://github.com/caslabuiowa/CBP-Collocation.

\subsubsection{Example 1 (ODE)}
Consider the ODE initial value problem (Lane-Emden equation of index 5)
\begin{equation} \label{eq:example1}
    \begin{split}
        x''(s)+\frac{2}{s}x'(s)+(x(s))^5 = 0, \\  
        x(0) = 1, \ x'(0) = 0,
    \end{split}
\end{equation}
over the interval $s\in[0,3]$, for which the exact solution is 
\begin{equation*}
    x(s) = \frac{1}{\sqrt{1+s^2/3}}.
\end{equation*}

\begin{figure}
    \centering
    \includegraphics[width=0.9\linewidth]{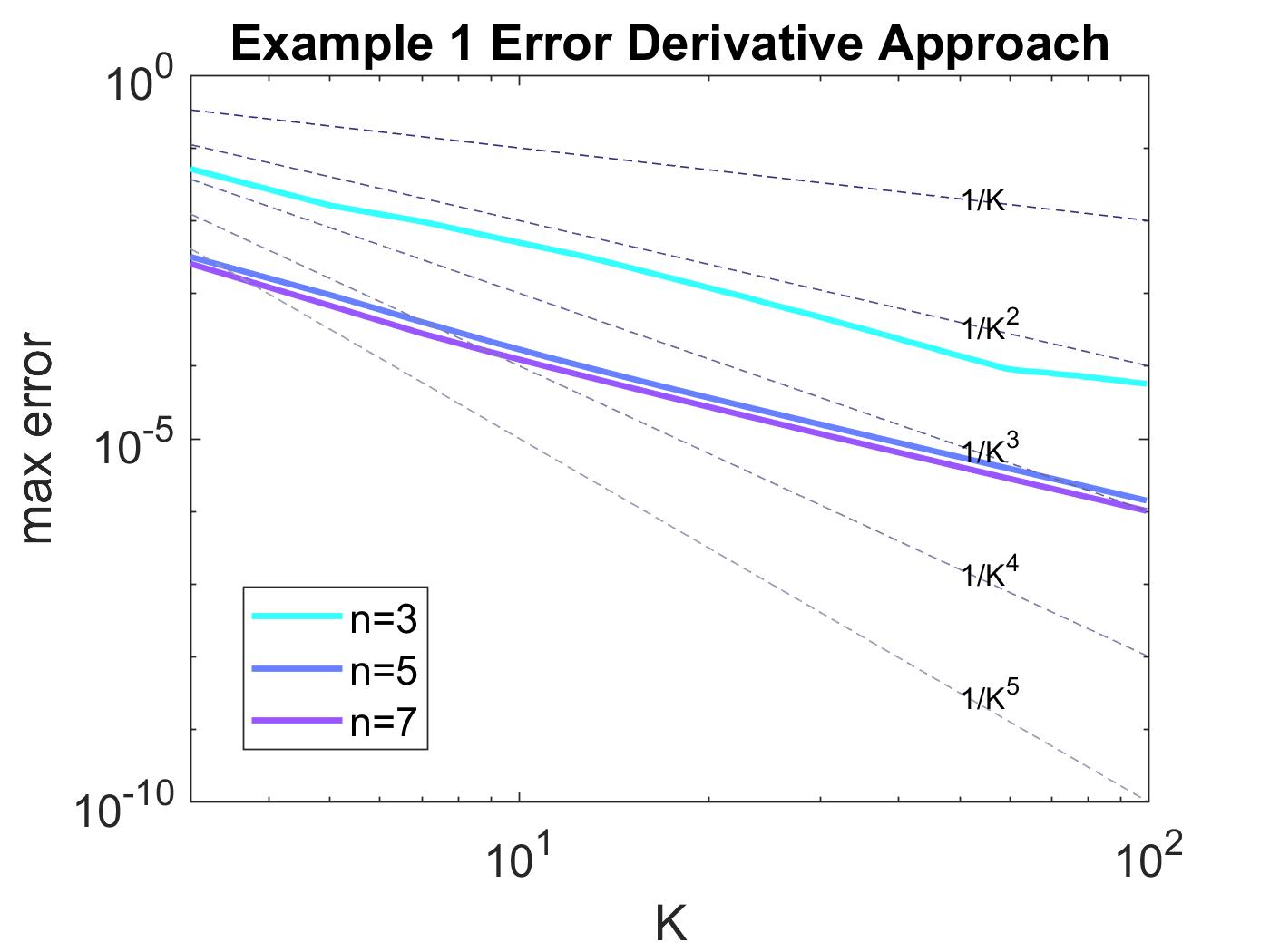}
    \caption{ Absolute error for derivative based approach to example 1.}
    \label{fig:ex1_der}
\end{figure}

First, we reformulate the problem to take the form of Problem \ref{prob:deriv}. The vector of unknowns in this case is $\bar{X}^{n,K}$, so we transform \eqref{eq:example1} to be discretized in terms of this decision variable 
\begin{equation}
\begin{split}
    \begin{split}
        \bar{X}^{n,K}(\textbf{D}_{n,K})^2+2\oslash s_{n,K}\circ(\bar{X}^{n,K}\textbf{D}_{n,K})+(\bar{X}^{n,K})^{\circ 5}=0, \\
        \bar{X}^{n[0]}_0=1, \\
        \{\bar{X}^{n[0]}\textbf{D}_{n}\}_0=0, \\
        \bar{X}^{n[j]}_n=\bar{X}^{n[j+1]}_0, \\
        \{\bar{X}^{n[j]}\textbf{D}_n\}_f=\{\bar{X}^{n[j+1]}\textbf{D}_n\}_0,\\
        j=1,\dots,K-2,
    \end{split}
\end{split}
\end{equation}
where \(\circ\) and \(\oslash\) denote the Hadamard product and division respectively. In the exponent, \(\circ\) signifies the use of Hadamard power. Figure \ref{fig:ex1_der} shows the error convergence rate of this method for increasing K and fixed \(n\).

\begin{figure}
    \centering
    \includegraphics[width=0.9\linewidth]{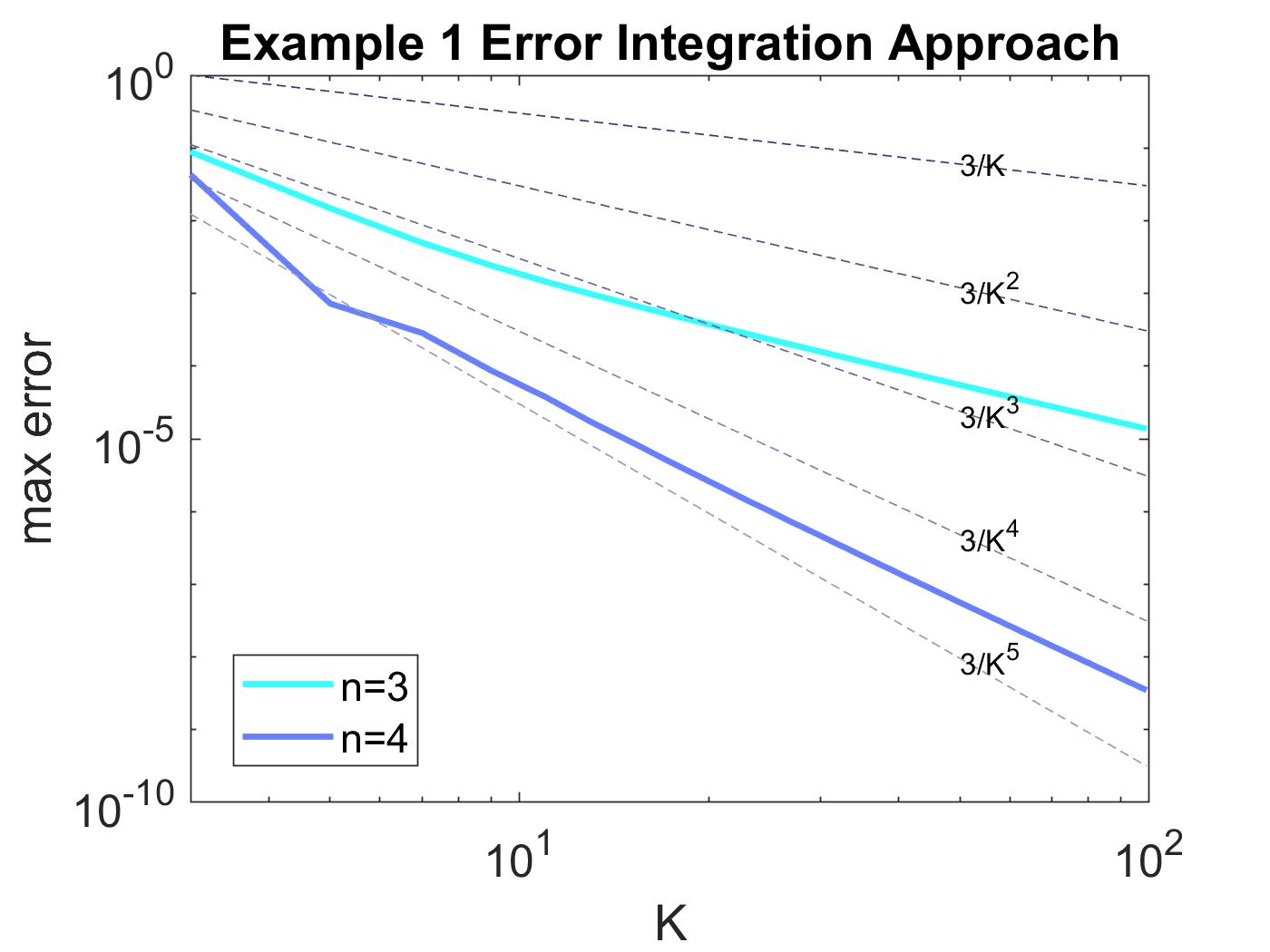}
    \caption{Absolute error for integration based approach to Example 1.}
    \label{fig:ex1_int}
\end{figure}

\begin{figure}
    \centering
    \includegraphics[width=0.9\linewidth]{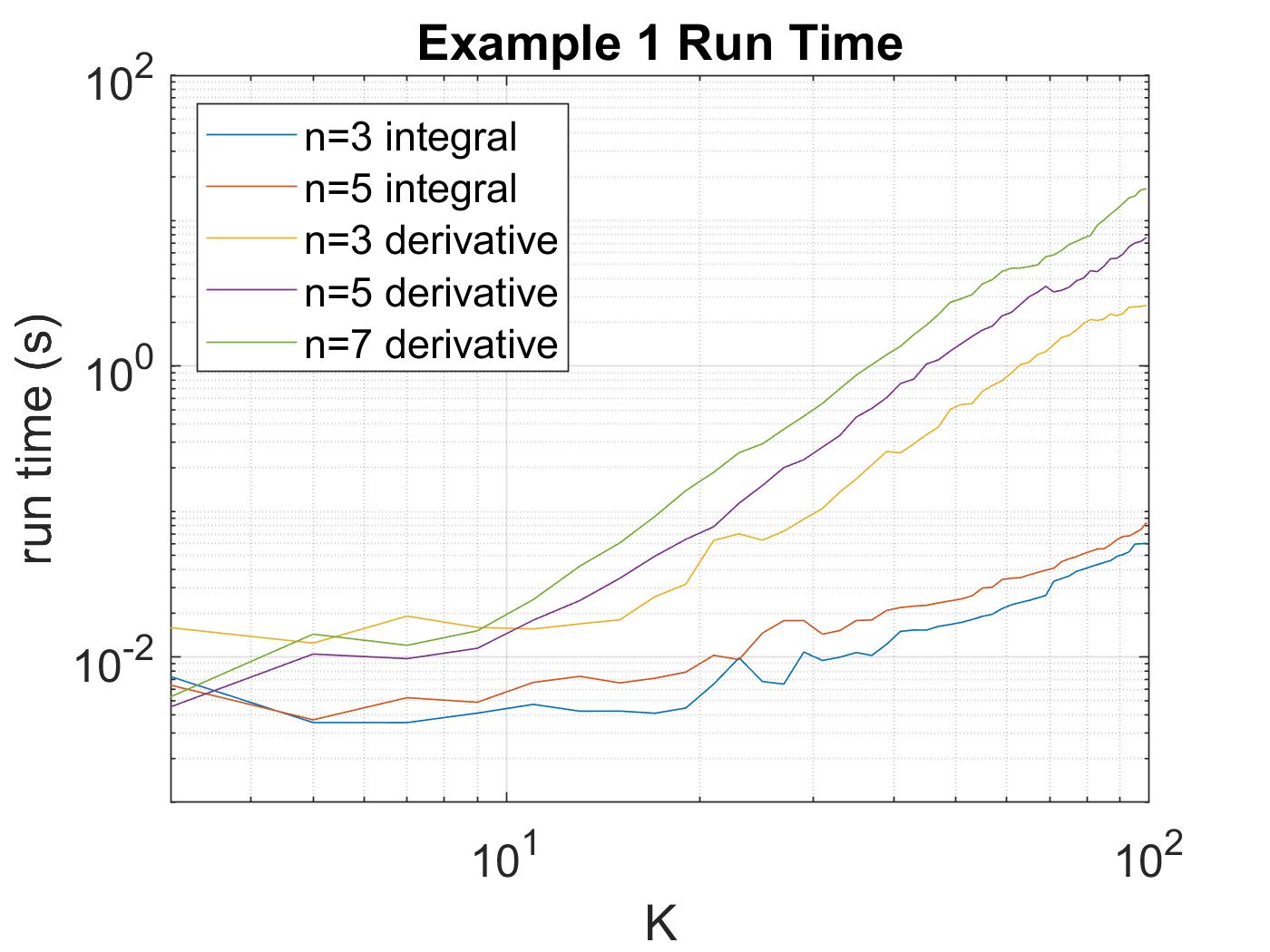}
    \caption{Runtime comparison between control point and knot only collocation approaches to Example 1.}
    \label{fig:ex1_run}
\end{figure}

Reformulating the problem to the form of Problem \ref{prob:integ} with $M=3$, the decision variable becomes $\theta_0=[\bar{X}^{(3),0,K}, X_0,X'_0]$. The problem can be written in terms of this variable as
\begin{equation}
    \begin{split}
        \begin{split}
            \theta_0\textbf{T}_2+2\oslash s_i\circ(\theta_0\textbf{T}_1)+(\theta_0\textbf{T}_0)^{\circ 5}=0, \\
            \{\theta_0\textbf{T}_0\}_0=1,\\
            \{\theta_0\textbf{T}_1\}_0=0.
        \end{split}
    \end{split}
\end{equation}
For $M=4$, we add the additional condition discussed in Remark \ref{rem:add_cond}
\begin{equation}
    \begin{split}
        X''(\sigma)+\frac{2}{\sigma}X'(\sigma)+(X(\sigma))^5=0 \, , \quad \sigma = \frac{s_1-s_0}{2} .
    \end{split}
\end{equation}
Figure \ref{fig:ex1_int} shows the error convergence rate of this method for increasing K.
\begin{remark}
    A substantial benefit of the integration approach is the reduction in runtime which accompanies its implementation as compared to the derivative approach. Figure \ref{fig:ex1_run} shows the difference in runtime for solutions of the same order to highlight this.  This increase in speed likely stems from a combination of the elimination of continuity constraints and a reduction in the decision variable size. It is also notable in this regard that the integration matrix is triangular and well conditioned compared to the derivative matrix.
\end{remark}

\vspace{6pt}
\subsubsection{Example 2 (ODE)}
\begin{figure}
    \centering
    \includegraphics[width=0.9\linewidth]{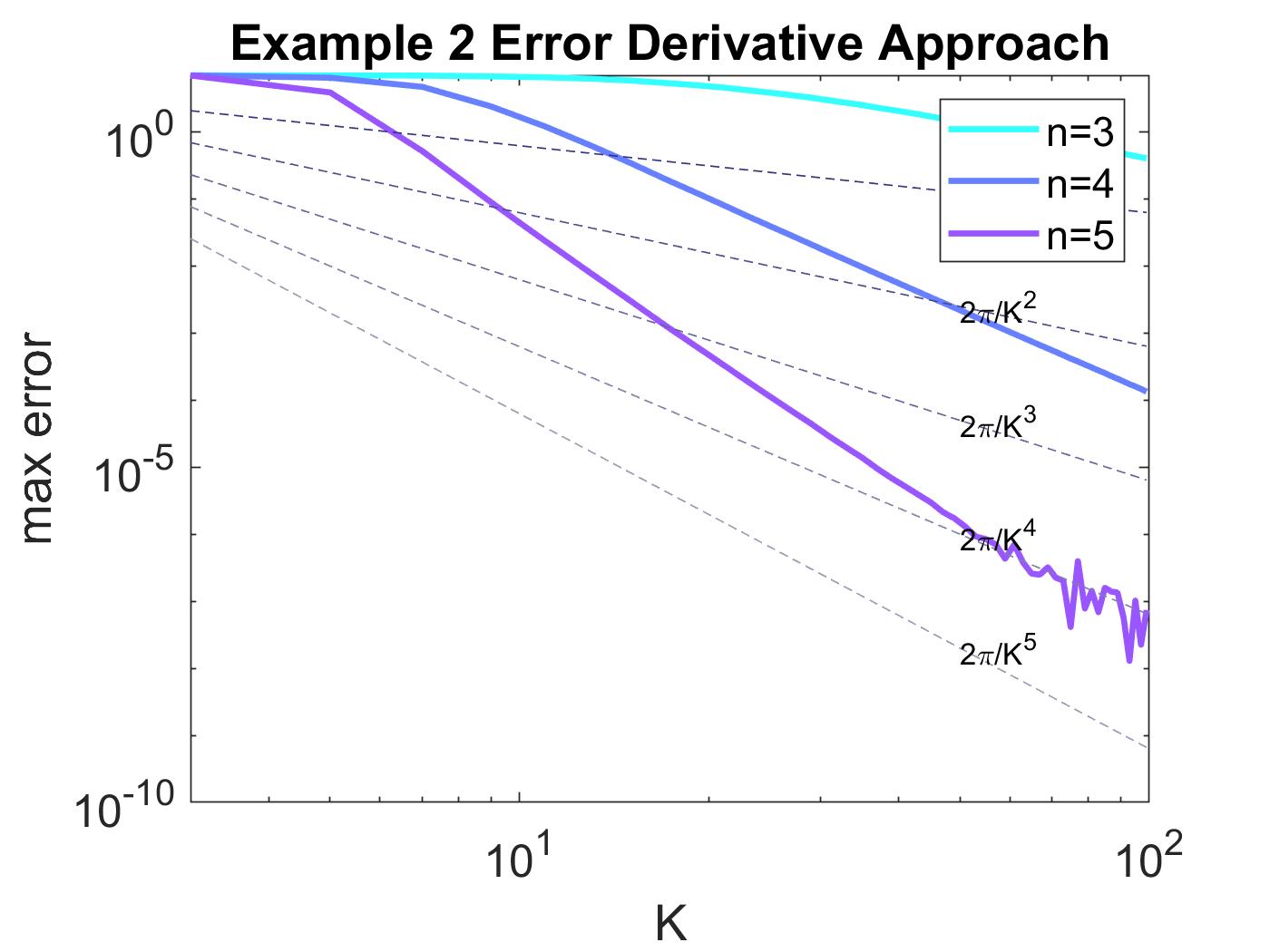}
    \caption{Absolute error for derivative based approach to Example 2. Note: Large oscillations here are a product of fsolve function tolerance options.}
    \label{fig:ex2_der}
\end{figure}

\begin{figure}
    \centering
    \includegraphics[width=0.9\linewidth]{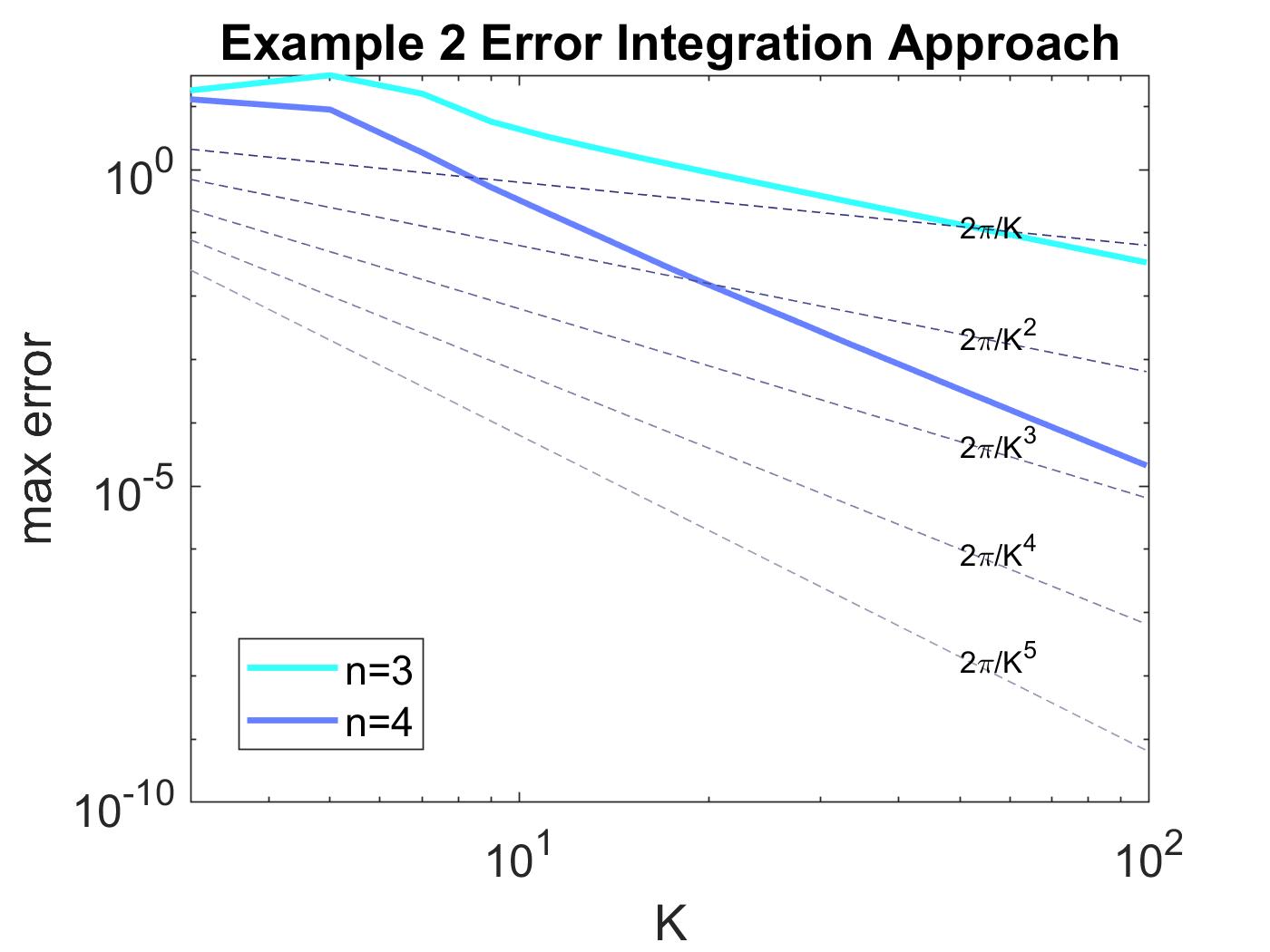}
    \caption{Absolute error for integration based approach to Example 2.}
    \label{fig:ex2_int}
\end{figure}

Consider the second order boundary value problem
\begin{equation}
    x''(s)+3x(s)=0, \ x(0)=7, \ x(2\pi)=0,
\end{equation}
on the interval $s\in[0,2\pi]$, for which the solution is 
\begin{equation*}
    x(s)=7\cos(\sqrt{3}s)-7\cot(2\sqrt{3}\pi)\sin(\sqrt{3}s).
\end{equation*}
Error for the derivative based implementation is given in Figure \ref{fig:ex2_der} and integration based implementation is given in Figure \ref{fig:ex2_int}.

\vspace{6pt}
\subsubsection{Example 3 (ODE)}
\begin{figure}
    \centering
    \includegraphics[width=0.9\linewidth]{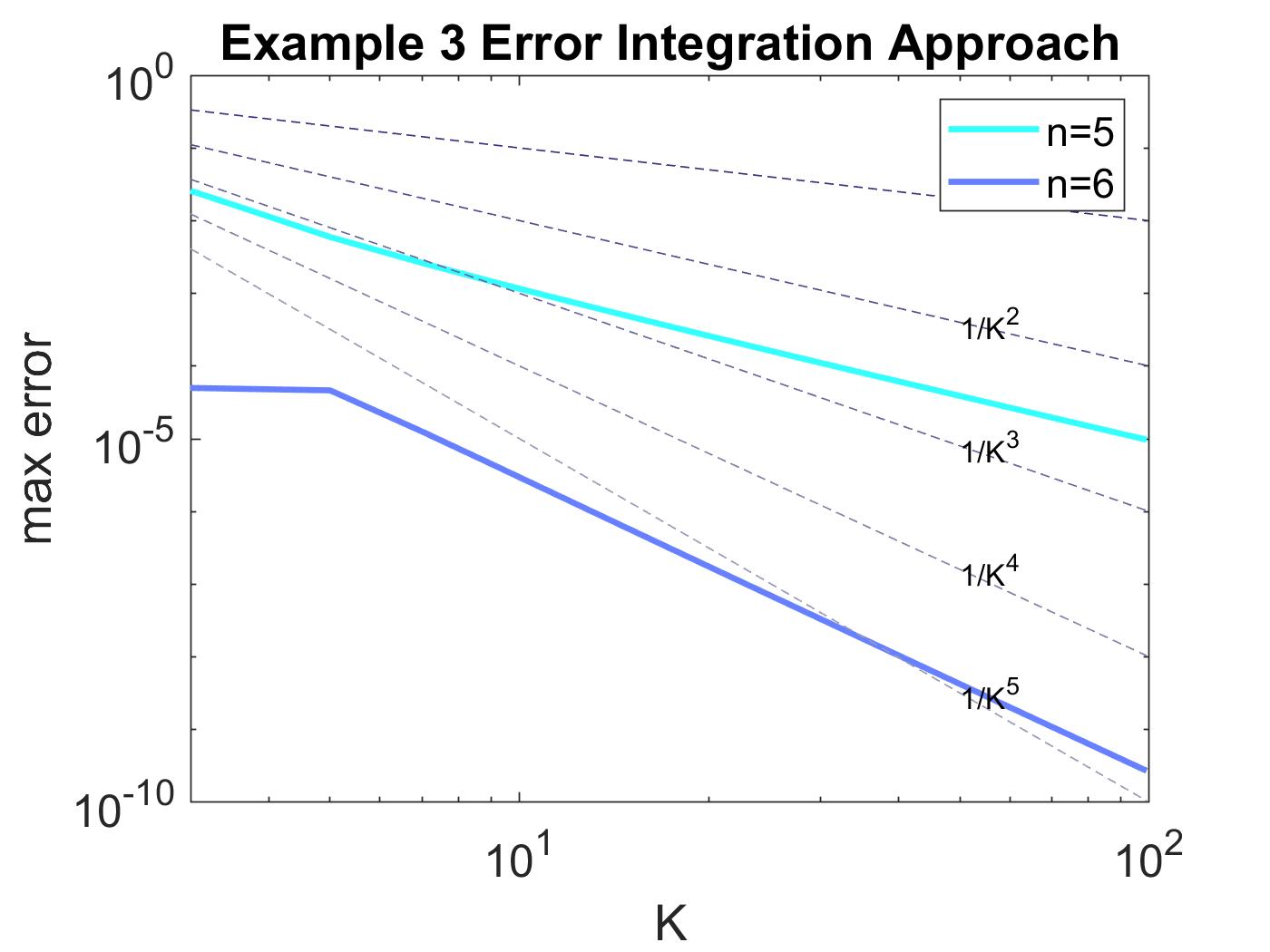}
    \caption{Absolute error for integration based approach to Example 3.}
    \label{fig:ex3_int}
\end{figure}

Consider the higher-order linear ODE
\begin{equation}
    \begin{split}
        x^{(4)}(s)-5x''(s)+4x(s)=\sin(s)+\cos(2s), \\
        x(0)=-1, \ x'(0)=0, \ x''(0)=-2, \ x^{(3)}=1,
    \end{split}
\end{equation}
on the interval $s\in[0,1]$, for which the exact solution is
\begin{equation*}
\begin{split}
        x(s) = \frac{1}{240}(61e^{-2s}-52e^{-s}-92e^s+29e^{2s} \\
        -6\sin^2(s)+24\sin(s)+6\cos^2(s)).
\end{split}
\end{equation*}
Error for the integration based implementation is given in Figure \ref{fig:ex3_int}.

\vspace{6pt}
\subsubsection{Example 4 (ODE)}
\begin{figure}
    \centering
    \includegraphics[width=0.9\linewidth]{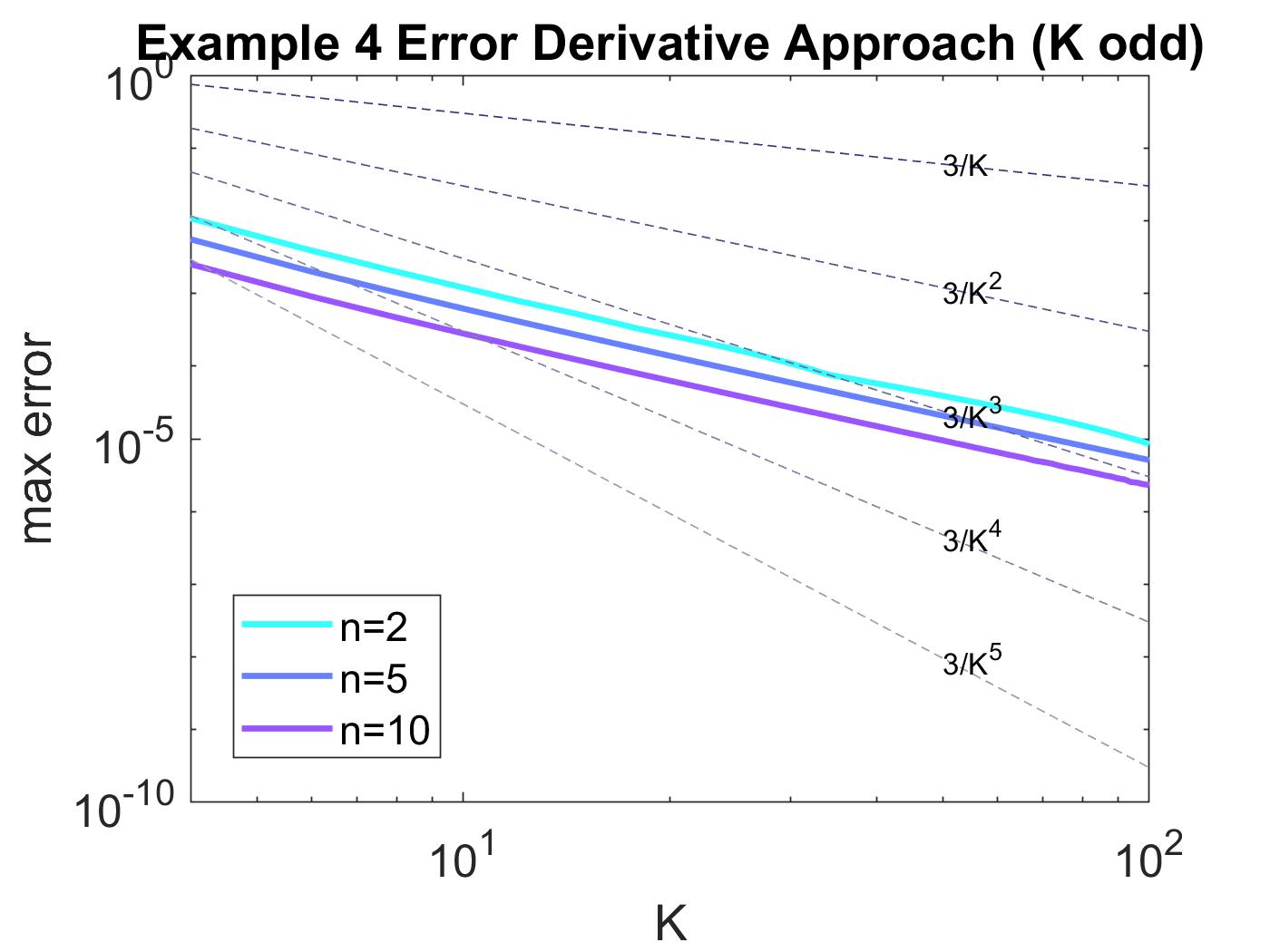}
    \caption{Absolute error for derivative based approach to Example 4 with odd $K$.}
    \label{fig:ex4_der_odd}
\end{figure}

\begin{figure}
    \centering
    \includegraphics[width=0.9\linewidth]{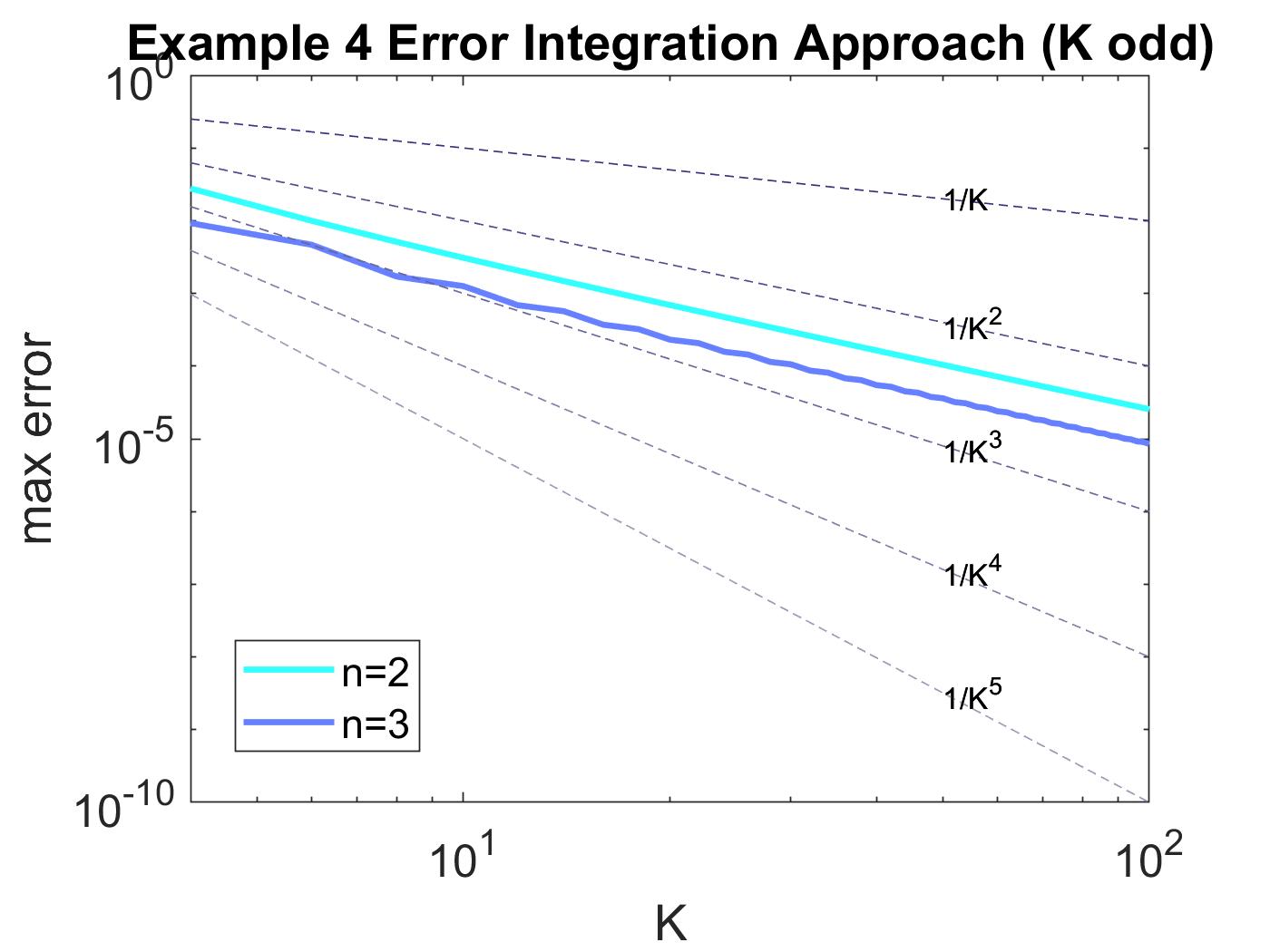}
    \caption{Absolute error for integration based approach to Example 4 with odd $K$.}
    \label{fig:ex4_int_odd}
\end{figure}

\begin{figure}
    \centering
    \includegraphics[width=0.9\linewidth]{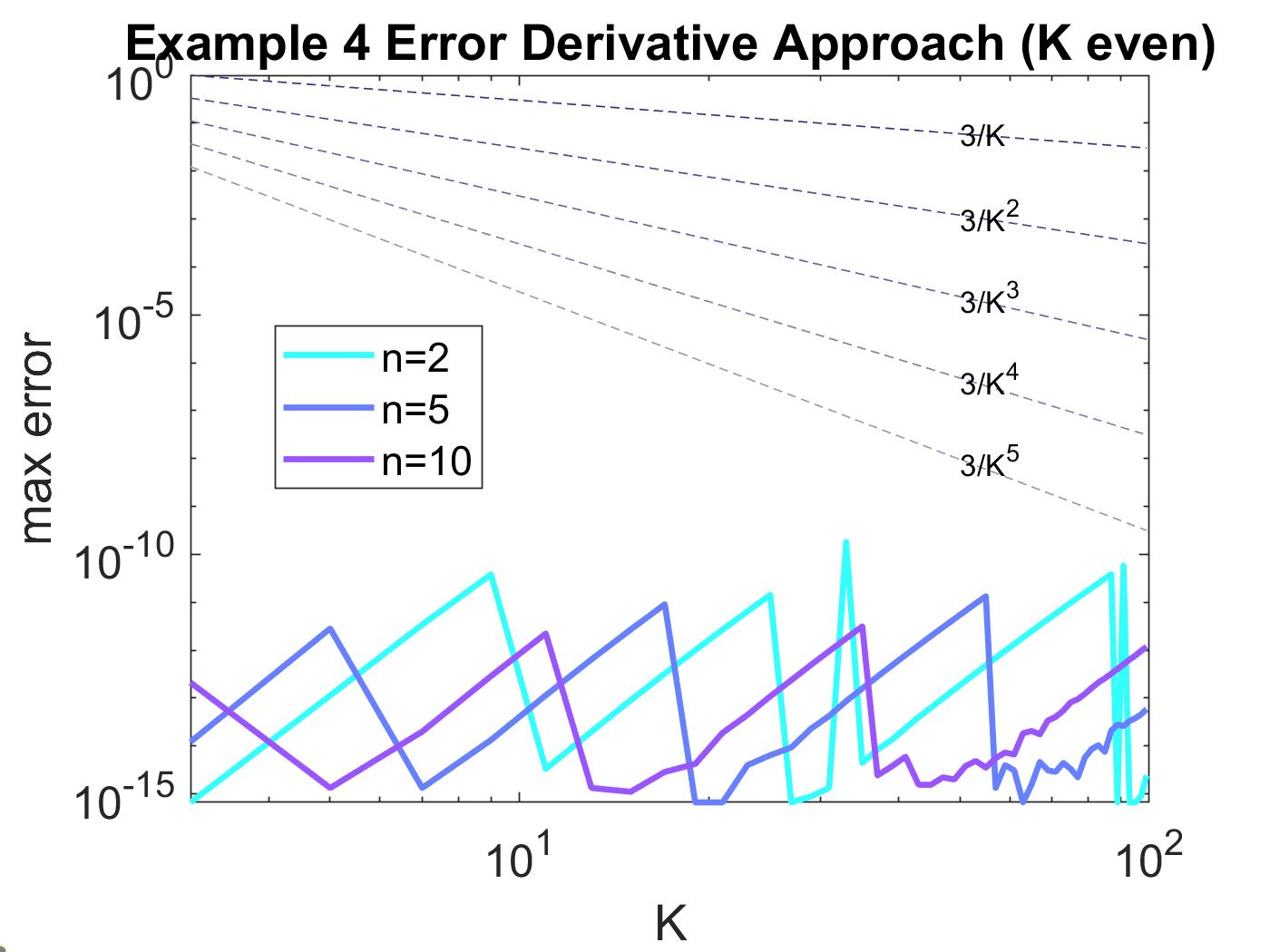}
    \caption{Absolute error for derivative based approach to Example 4 with even $K$.}
    \label{fig:ex4_der_even}
\end{figure}

\begin{figure}
    \centering
    \includegraphics[width=0.9\linewidth]{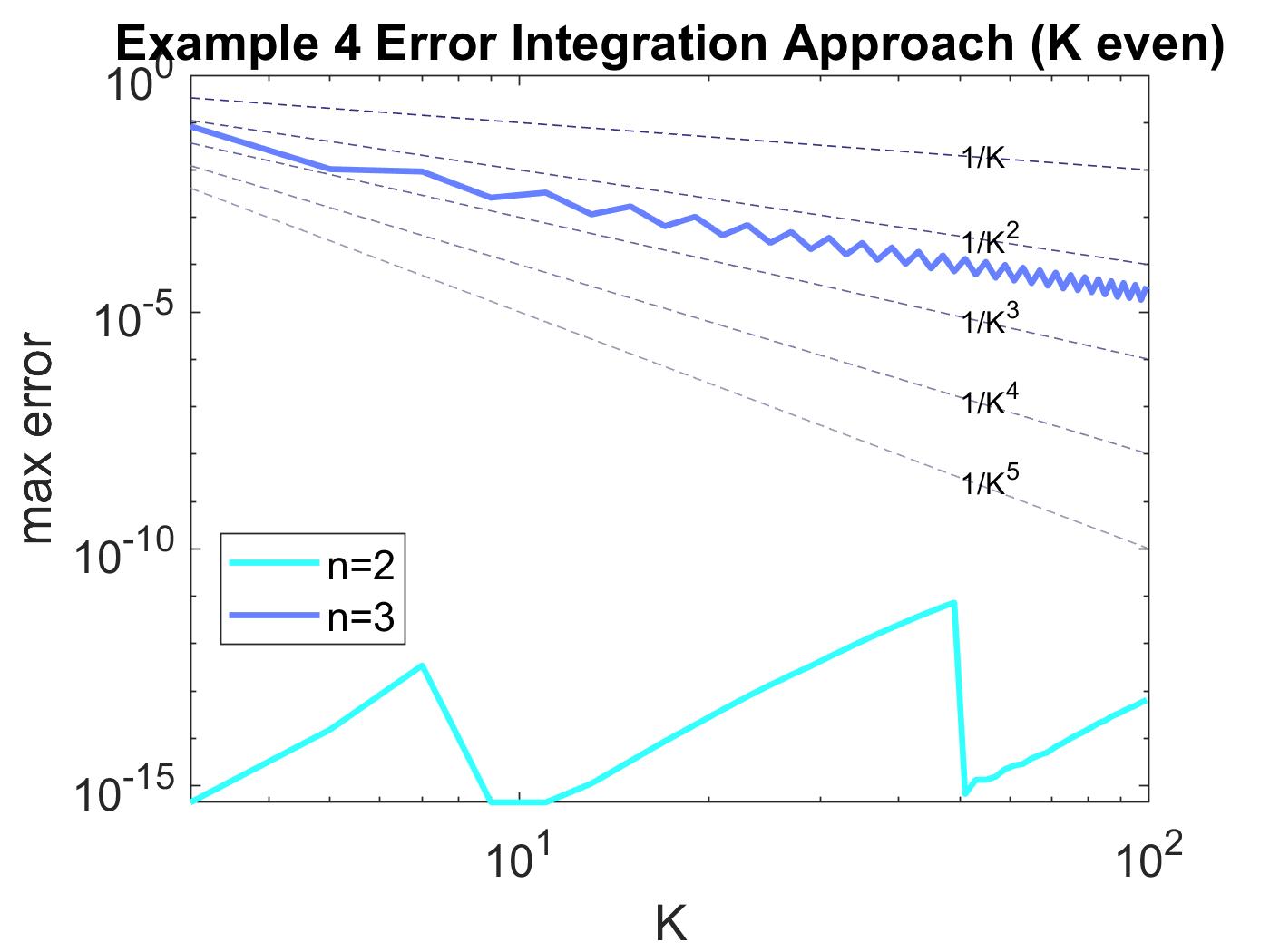}
    \caption{Absolute error for integration based approach to Example 4 with even $K$.}
    \label{fig:ex4_int_even}
\end{figure}

Consider the ODE
\begin{equation}
    x'-|s-0.5|=0, \ x(0)=1,
\end{equation}
over the interval $s\in[0,1]$, for which the $C^1$ continuous solution is 
\begin{equation}
    \begin{cases}
        x(s)=-\frac{1}{2}s^2+\frac{1}{2}s+1, \ s\in[0,0.5]\\
        x(s)=\frac{1}{2}(\frac{1}{2}-s)^2+1.125, \ s\in[0.5,1]
    \end{cases}.
\end{equation}
Note, the discontinuity on $x''(s)$ will occur at $s=0.5$. For equidistant knot placement on the given interval, this will lead to an exact solution to the problem for even $K$ and a convergence rate similar to other problems for odd $K$ (Figure \ref{fig:ex4_der_odd} and Figure \ref{fig:ex4_int_odd}). For the derivative case, the solution will remain exact for increasing $n$ (Figure \ref{fig:ex4_der_even}). However, the integration method will only obtain an exact solution for $M=2$ as a result of implicit continuity enforcement  (Figure \ref{fig:ex4_int_even}).

\vspace{6pt}
\subsubsection{Example 5 (DAI)}
Consider the problem of a dynamic system maintaining a bounded distance from a moving target
\begin{equation}
    \begin{split}
        \ddot{x}(t) = u(t)-\gamma\dot{x},\\  
        x(0) = 0, \ \dot{x}(0) = 0,\\
        |u(t)| \le u_{max}, \\
        |x(t)-x_{des}(t)| \le r, \\ 
        t\in[0,T].
    \end{split}
\end{equation}
Here, let the target function $x_{des} = 2\sin(t)\cos(2t)$, the system's damping coefficient $\gamma = 1$, the upper limit on the input $u_{max} = 5$, and the maximum distance from the target function $r=0.2$.

\begin{figure}
    \centering
    \includegraphics[width=\linewidth]{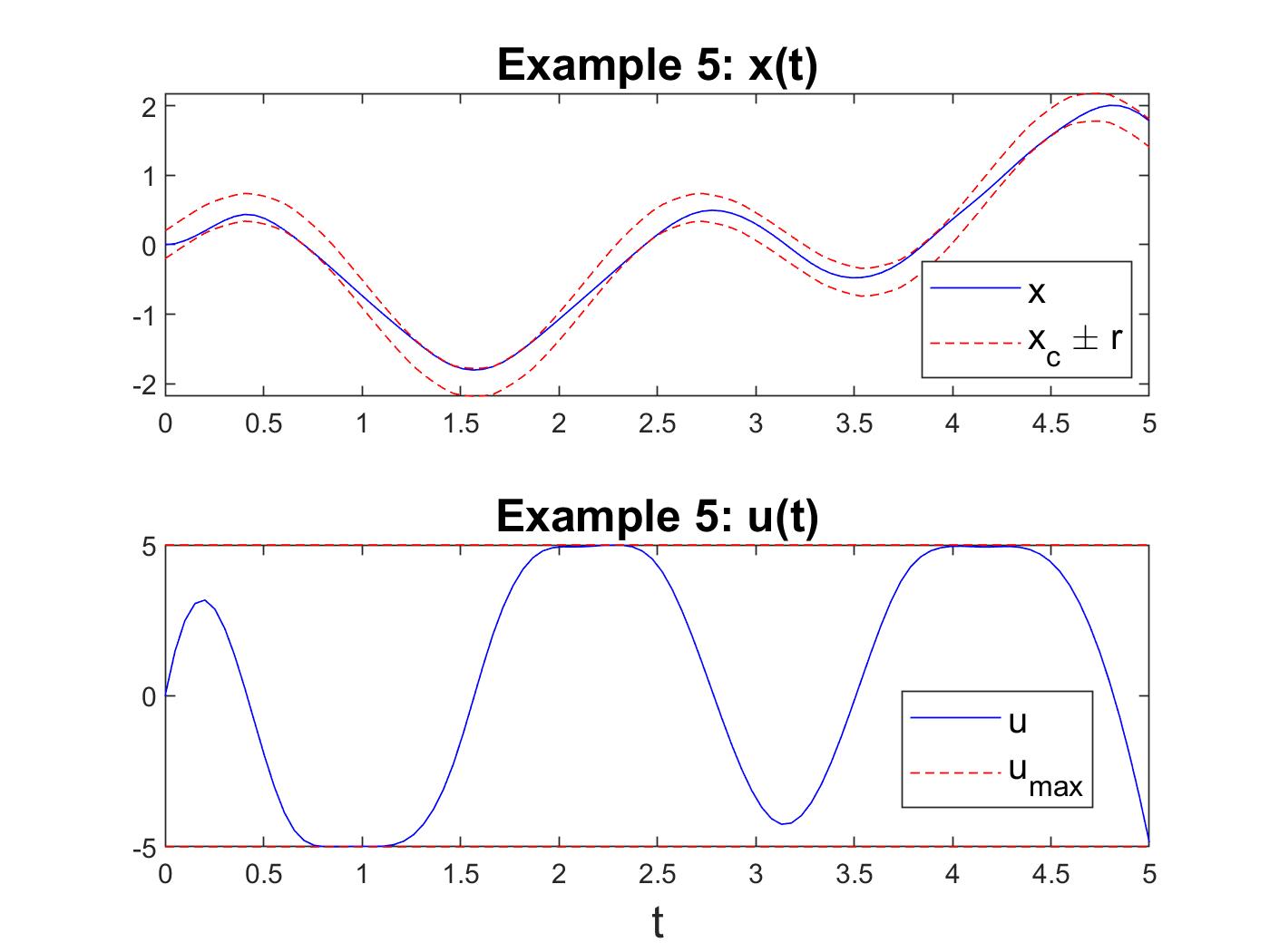}
    \caption{Example solution to Example 5 with $K=25$.}
    \label{fig:ex5_int}
\end{figure}

Employing the integration solution method, let $M = 3$ and the decision variable be $[\theta_{0,x},\theta_{0,u}]$ representing the variables $x$ and $u$ respectively. The problem can be written in discrete from in terms of equality and inequality constraints. From the dynamics we have the equalities
\begin{equation*}
    \begin{split}
        \theta_{0,x} \textbf{T}_2 = \theta_{0,u} \textbf{T}_0 - \gamma \theta_{0,x} \textbf{T}_1 \\
        \{\theta_{0,x}\textbf{T}_0\}_0=0,\\
        \{\theta_{0,x}\textbf{T}_1\}_0=0.    
    \end{split}
\end{equation*}
Notably, the equality constraints are only enforced at the knots with this methodology; however, the inequality constraints must be checked on the full vector of control points. For this problem, the inequality constraints can be written
\begin{equation*}
    \begin{split}
        |\bar{U}^{n,K}_{i,j}|-u_{max} \le 0, \\
        |\bar{X}^{n,K}_{i,j}-\bar{X}^{n,K}_{des,i,j}| - r \le 0,
    \end{split}
\end{equation*}
where $\bar{U}^{n,K}, \ \bar{X}^{n,K}$ and $\bar{X}^{n,K}_{des}$ are the control points of the CBPs in the form of equation \eqref{eq:comp} approximating $u(t), \ x(t)$ and $x_{des}(t)$ respectively. These constraints can be given to an off the shelf optimizer for a solution and an example of one such solution is given in Figure \ref{fig:ex5_int}.

\vspace{6pt}
\subsubsection{Example 6 (DAI)}
\begin{figure}
    \centering
    \includegraphics[width=\linewidth]{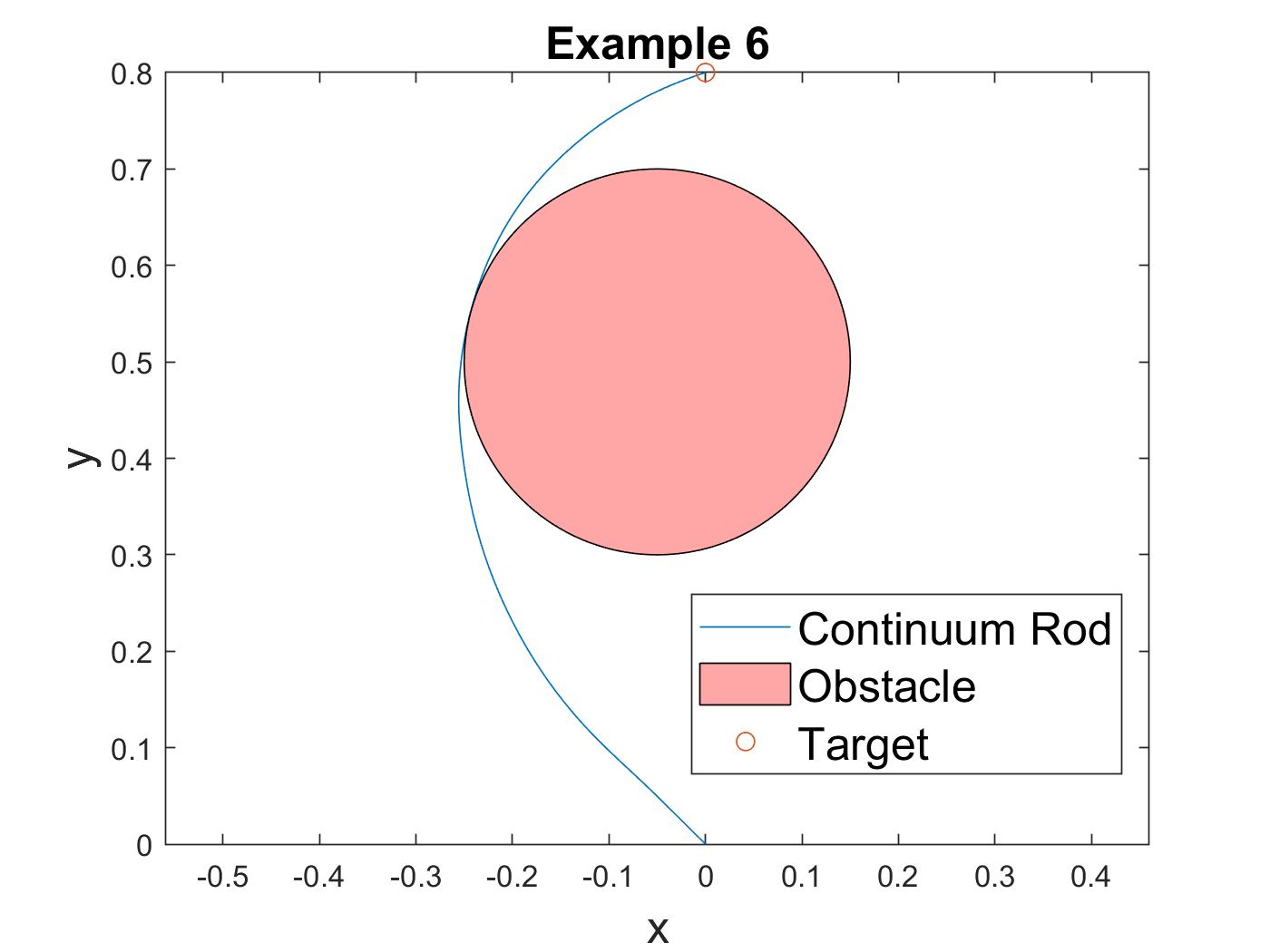}
    \caption{Example solution to Example 6 with $K=15$.}
    \label{fig:ex6_int}
\end{figure}

Consider the problem of a 2D continuum rod reaching around an object to a desired end point $p_{des}$
\begin{equation}
    \begin{split}
        p'(s) = R(s)\nu(s), \\
        R'(s) = R(s)u(s), \\
        p(0) = p_0,\\
        p(L) = p_{des}, \\
        ||p(s)-O||_2 \ge r_{safe}, \\
        \nu_{min} \le \nu(s) \le \nu_{max}, \\
        |u(s)| \le u_{max}.
    \end{split}
\end{equation}
In this problem, $p(s)$ and $R(s)$ represent the position and orientation of the rod respectively, and $\nu(s)$ and $u(s)$ are the corresponding linear and angular strains. $O$ is the origin of a circular object, and $r_{safe}$ is the object's radius. In 2D, the orientation $R(s)$ can be parameterized by a single angle $\phi(s)$ such that
\begin{equation*}
    R(\phi(s)) = \begin{bmatrix}
        \cos(\phi(s)) & -\sin(\phi(s)) \\
        \sin(\phi(s)) & \cos(\phi(s))
    \end{bmatrix} .
\end{equation*}
Approaching with the integration method, the decision variable is chosen as  $[\theta_{0,\nu},\theta_{0,\phi}]$ leading to
\begin{equation*}
    \begin{split}
        \bar{P}_0^{n+1,[0]} = p_0, \ \bar{P}_{n+1}^{n+1,[K]} = p_{des}, \\
        \nu_{min} \le \bar{V}^{n,K}_{i,j} \le \nu_{min}, \\
        |\bar{U}^{n,K}_{i,j}| \le u_{max}, \\
        ||\bar{P}_{i,j}^{n+1,K}-O||_2^2 \ge r_{safe}^2,
    \end{split}
\end{equation*}
where $\bar{P}^{n+1,K}$ are the control points of the CBP approximation of $p(s)$ computed from the decision variable, and $\bar{V}^{n,K}$ and $\bar{U}^{n,K}$ are the control points of the linear and angular strain functions in form \eqref{eq:comp}. An example solution can be seen in Figure \ref{fig:ex6_int}.

\section{conclusion}  \label{sect:conc}
In this paper we introduce two collocation methods with unique benefits for solving ODEs which yield solutions in the form of composite Bernstein polynomials  as well as an extension for the solution of DAIs. This polynomial basis has desirable properties for a number of fields, particularly in optimal control where these methods can be applied for the efficient approximation of system dynamics and constraints within OCPs. In fact, solutions to DAIs as provided here may serve as strong initial guesses in OCPs with the same dynamics and constraints. We provide a number of numerical examples of both ODEs and DAIs to demonstrate the effectiveness of the proposed methodology.

\appendices
\section{} \label{app:int}
The integration matrix used in \eqref{eq:theta_int} is given by
\begin{equation}
    \boldsymbol{\zeta}_{n}=\begin{bmatrix}
         \boldsymbol{\Gamma}_{n} & \textbf{0}\\
        \boldsymbol{\Psi}_n & \textbf{I}_M
    \end{bmatrix}\\    
\end{equation}
where $\boldsymbol{\zeta}_{n}\in\mathbb{R}^{K(n+1)+M\times K(n+2)+M},$ and
\begin{equation*}
    \begin{split}
        \boldsymbol{\Gamma}_{n} = \begin{bmatrix}
                \boldsymbol{\gamma}_n^{[0]} & \boldsymbol{\phi}_{0,n} & \dots & \dots &  \boldsymbol{\phi}_{0,n} \\
                0 & \boldsymbol{\gamma}_n^{[1]} & \boldsymbol{\phi}_{1,n} & \dots & \boldsymbol{\phi}_{1,n} \\
                \vdots & \ddots & \ddots  & \ddots & \vdots \\
                \vdots & \ddots & \ddots & \ddots & \boldsymbol{\phi}_{K-2,n} \\
                0 & \dots & \dots & 0 & \boldsymbol{\gamma}_n^{[k-1]}
        \end{bmatrix},  \\
        \boldsymbol{\phi}_{i,n} = \frac{s_{i+1}-s_i}{n+1}\begin{bmatrix}
            1 & \dots & 1 \\
            \vdots & \ddots & \vdots \\
            1 & \dots & 1
        \end{bmatrix}\in \mathbb{R}^{n+1\times n+2}, \\
        \boldsymbol{\Psi}_n = \{\psi_{i,j}\} \in\mathbb{R}^{M\times K(n+2)},\\
        \psi_{i,j} = \begin{cases}
            1 \quad i = n+1 \\
            0 \quad \text{otherwise}
        \end{cases},
    \end{split}
\end{equation*}
$\boldsymbol{\gamma}_n^{[i]}$ is the single integration matrix for $[s_i,s_{i+1}]$ introduced in Equation \eqref{eq:single_int} and reported below for convenience, 
\begin{equation*} 
        \boldsymbol{\gamma}_n = \frac{s_f-s_0}{n+1}\begin{bmatrix}
        0 & 1 & \dots & 1 \\
         \vdots  & \ddots & \ddots &\vdots \\
        0  & \dots  & 0 & 1
        \end{bmatrix}  \in \mathbb{R}^{n+1\times n+2},
\end{equation*}
$\textbf{I}_M$ is the identity matrix in $\mathbb{R}^{M\times M}$, and $\boldsymbol{\Gamma}_{n}\in\mathbb{R}^{K(n+1)\times K(n+2)}$. Importantly, the structure of $\boldsymbol{\Gamma}_{n}$ inherently enforces continuity at the knots of the resulting CBP. 

Next, we show how matrix \(\boldsymbol{\zeta}_{n}\) is derived. Consider the vector of control points for the $m$-th derivative of some composite Bezier polynomial (CBP), denoted as \(X^{(m)}(s)\), specifically, \(\bar{X}^{(m),n,K}\). We will focus on the \(i\)th BP within the CBP, represented by \(x_n^{(m),[i]}(s)\) for \(s \in [s_i, s_{i+1}]\), and its corresponding vector of control points \(\bar{x}^{(m),n[i]}\). Utilizing Property \ref{prop:int}, we determine the control points of the integral, \(x_{n+1}^{(m-1),[i]}(s)\), represented as \(\bar{x}^{(m-1),n+1 [i]}\). These are computed using the formula:
\begin{equation} \label{eq:int_polyi}
    \bar{x}^{(m-1),n+1 [i]} = \bar{x}^{(m),n[i]} \boldsymbol{\gamma}_n^{[i]} + \bar{x}_0^{(m-1),[i]} \mathbbm{1}_{n+2}^\top.
\end{equation}
To extend this from the individual polynomial case to the composite case, recall the structure of the vector $\theta_n= [\bar{X}^{(m),n,K},X_0,\dots,X^{(m-1)}_0]$. In the operation \(\theta_{n+1}=\theta_{n}\boldsymbol{\zeta}_{n}\), the matrix \(\boldsymbol{\zeta}_{n}\) adds the the appropriate initial condition to all control points in \(\bar{X}^{(m-1),n+1,K}\) with the submatrix \(\boldsymbol{\Psi}_n\). The initial conditions are persevered in the resulting vector \(\theta_{n+1}\) with the identity matrix \(\textbf{I}_m\). The first term of \eqref{eq:int_polyi} is captured by \(\boldsymbol{\Gamma}_n\) with the addition of a $C^0$ continuity constraint embedded in the matrix. Note that \(\bar{x}_0^{(m-1),[i]}\neq X^{(m-1)}_0\). For continuous functions, the upper triangle of \(\boldsymbol{\Gamma}_n\) is populated by \(\boldsymbol{\phi}_{i,n}\) which shifts the polynomial and ensures the equality \(\bar{x}^{(m-1),n+1 [i]}_{n+1}=\bar{x}^{(m-1),n+1 [i+1]}_0\), resulting in $C^0$ continuity through Property \ref{prop:end}.

The matrix $\textbf{P}_n$ used to extract the knots of a CBP from $\theta_n$ in \eqref{eq:transformation} is obtained by
\begin{equation}
    \begin{split}
        \textbf{P}_{n}=\{p_{ij}\}\in\mathbb{R}^{K (n+1)+M\times K+1},\\
        p_{i,j} = \begin{cases}
            1 \quad i = j = 0 \\
            1 \quad i = k(n+1)-1, \ j = k \\
            0 \quad \text{otherwise}
        \end{cases}, \\
        k = 1,\dots,K.
    \end{split}
\end{equation}
The following link provides the Matlab function that is used to compute the matrix $\boldsymbol{\zeta}_n$ in equation (42): https://github.com/caslabuiowa/CBP-Collocation. 

\bibliographystyle{ieeetr}
\bibliography{refs}
\end{document}